\documentstyle[12pt]{article}
\newtheorem{theor}{Theorem}

\newtheorem{prop}{Proposition}

\author{M.Deza, ENS, Paris, France \and V.P.Grishukhin, CEMI RAN, Moscow, Russia}

\title{Once more about the 52 four-dimensional parallelotopes}

\date{}

\begin{document}
\maketitle

\begin{abstract}
There are several works \cite{De} (and \cite{St}), \cite{En}, \cite{Co}
and \cite{Va} enumerating four-dimensional parallelotopes. In this work
we give a new enumeration showing that any four-dimensional parallelotope
is either a zonotope or the Minkowski sum of a zonotope with the regular
24-cell $\{3,4,3\}$. Each zonotopal parallelotope is the Minkowski sum
of segments whose generating vectors form a unimodular system. There
are exactly 17 four-dimensional unimodular systems. Hence, there are 17
four-dimensional zonotopal parallelotopes. Other 35 four-dimensional
parallelotopes are: the regular 24-cell $\{3,4,3\}$ and 34 sums of the
regular parallelotope with non-zero zonotopal parallelotopes. For the
nontrivial enumerating of the 34 sums we use a theorem describing
necessary and sufficient conditions when the Minkowski sum of a
parallelotope with a segment is a parallelotope.
\end{abstract}

\section{Introduction}
A {\em parallelotope} is a convex polytope which fills the space facet to
facet by its translation copies without intersecting by inner points.
Such a filling by parallelotopes is a {\em tiling}. The centers of the
tiles form a lattice. A parallelotope of dimension $n$ is called
{\em primitive} if exactly $n+1$ adjacent parallelotopes meet in each
vertex of its lattice tiling. Voronoi defines {\em an L-type} of a
parallelotope, which is (in modern terms) the isomorphism class of the
face lattice of the parallelotope. A special kind of a parallelotope is
the Voronoi polytope of a lattice. The Voronoi polytope at a lattice
point $v$ is the set of points which are at least as close to $v$ as
to any other lattice point. Voronoi conjectured that every parallelotope
of a class of parallelotopes of the same L-type is affinely equivalent
to a Voronoi polytope (of course, of the same L-type), and
he proved this conjecture for primitive parallelotopes.

In four-dimensional space, Voronoi \cite{Vo} determined all the 3 types
of primitive parallelotopes. Using projection along a zone of parallel
edges, Delaunay \cite{De} find 51 types of four-dimensional
parallelotopes. The missed 
52th type was discovered by Shtogrin \cite{St}.

Engel verified by computer the result 
of Delaunay corrected by Shtogrin.
In Table 1 of \cite{En} he gives an informative 
and useful list of
parallelotopes of the 52 types, each of maximal symmetry. Besides, 
in
Fig.1 he gives a partial order between the 52 parallelotopes. This order
consists 
of two disjoint components. These two partial orders
correspond to partial orders 
between zonotopes: the first one amongst zonotopes, which are themselves 
parallelotopes, and the second one amongst
zonotopes, the Minkowski sum of which with the 24-cell gives a parallelotope.
Conway in the chapter "Afterthoughts: Feeling the form of a
Four-dimensional lattice" 
of \cite{Co} proposes {\em conorms} in order to characterize shapes (i.e. types) of 
parallelotopes. (In fact, Conway
enumerates shapes of four-dimensional Voronoi polytopes.) The 17 types of
parallelotopes which are zonotopes are parameterized by 16 subgraphs
of the complete graph $K_5$ and by the complete bipartite graph $K_{3,3}$.
The remaining 35 types of parallelotopes are characterized by shapes
of positions minimal conorms in $4\times 4$ matrices of all conorms.
Vallentin \cite{Va} repeated the above Conway's computations. He compute
in details the 35 nonequivalent conorms.

In \cite{Ry} Ryshkov asserts that a four-dimensional parallelotope is
either the Minkowski sum of segments, or the Minkowski sum of the
24-cell with a set of segments. The segments of the set are parallel
to the edges of the 24-cell. But no proof of these assertions was
published.

A zonotope is the Minkowski sum of segments. McMullen \cite{McM} proved
that an 
$n$-dimensional zonotope is a parallelotope if and only if the
directing vectors of 
its segments form an $n$-dimensional unimodular
system. It is well known (see, 
for example, \cite{DaG}) that there are
two maximal four-dimensional unimodular 
systems. These are the graphic
system of 10 vectors representing the regular 
matroid of the complete
graph $K_5$ and the cographic system of 9 vectors 
representing the
cographic matroid of the complete bipartite graph $K_{3,3}$. 
Each
subsystem of the maximal cographic unimodular system $K^*_{3,3}$ is
graphic 
and represents a subgraph of $K_5$. Hence, each four-dimensional
zonotopal 
parallelotope is generated either by one of 16 graphic unimodular
subsystems of 
$K_5$ or by the cographic system $K^*_{3,3}$.
Note that the graph $K_4$ in the table 
on page 87 of \cite{Co}
generates a three-dimensional unimodular system. But the rank 
4 subgraph
$C_{221}+\bf 1$ of $K_5$ is missed in the Conway's table. (An explanation

of the notation $C_{221}+\bf 1$ is given below in Section 3.)
Note that in \cite{Va},
 zonotopes are described by their facet vectors.
The facet vectors of zonotope $Z(M)$, related to a matroid $M$, form
a representation of the dual matroid $M^*$. Hence, in contrast to us,
Vallentin considers cographic matroids of corresponding graphs.

Remaining 35 parallelotopes are Minkowski sums of zonotopes related to
unimodular systems with the regular 24-cell $\{3,4,3\}$, and the 24-cell
itself. Note that zonotopes with the same unimodular system can give
distinct sums with the 24-cell.

In this paper we give a list of the 52 types of four-dimensional
parallelotopes, where for each type we explicitly give the corresponding
unimodular system. The case of four-dimensional parallelotopes is a very
instructive example of constructions of a large class of $n$-dimensional
parallelotopes.

\section{Parallelotopes of non-zero width}
Venkov introduced in \cite{Ve1} a notion of a polytope of non-zero
width in direction of a $k$-dimensional subspace $X^k$ as a polytope
whose intersection with any affine $k$-space parallel to $X^k$ is
either $k$-dimensional or empty. He studied parallelotopes of non-zero
width. The most interesting are parallelotopes of non-zero width in
direction of a line (or in direction of a vector spanning this line).

It is not difficult to see that if a parallelotope $P$ has a non-zero
width in direction of a line $l$, then the line $l$ is parallel to
some edges of $P$. A set of mutually parallel edges of $P$ is called
an {\em edge zone} of $P$. Following to Delaunay \cite{De}, Engel
called an edge zone {\em closed} if each two-dimensional face of $P$
has either two or none of edges of this zone. Otherwise, the zone is
called {\em open}.

For each edge zone $E$, there is a vector $z(E)$ with integer coordinates
which is parallel to the edges of $E$. Voronoi called the vector $z(E)$
a {\em characteristic} of any edge of the edge zone $E$. So, any edge
$e\in E$ has the form $e=\rho_ez(E)$. The positive number $\rho_e$ is
called by Voronoi {\em regulator} of the edge $e$.

In \cite{Gr}, the following proposition is proved.
\begin{prop}
\label{Pnz}
For a 
parallelotope $P$ and a vector $z$, the following assertions
are equivalent:

(i) $P$ has a closed edge zone parallel to $z$;

(ii) $P$ has a non-zero width in direction of $z$;

(iii) $P$ is the Minkowski sum of a segment $S(z)$ of the line spanned
by $z$ and a parallelotope $P'$ of zero width in direction $z$, i.e.,
$P=P'+S(z)$.
\end{prop}
We say that a parallelotope $P$ is {\em of zero width in direction}
$z$ if $P$ is not of non-zero width in this direction. The length of
the segment $S(z)$ in item (iii) is equal to the length of the shortest
edge of the closed zone parallel to $z$. It implies that this zone is
open in the parallelotope $P'$. Note that $S(z)=\lambda(z-z)$ for some
real $\lambda$, where $z-z=\{x=\alpha z: -1\le\alpha\le 1\}$ is the
Minkowski sum of $z$ and $-z$.

If the parallelotope $P'$ in the sum $P=P'+S(z)$ is also of non-zero
width in another direction $z_1$, then, by Proposition~\ref{Pnz}(iii),
$P'=P''+S(z_1)$. Of course, the Minkowski sum is associative and
distributive. Hence, if a parallelotope $P$ has a non-zero width in
several directions $z\in Q$, then $P=P_0+\sum_{z\in Q}S(z)$, where
$P_0$ is a parallelotope, which has no direction of non-zero width.
The sum $Z(Q)=\sum_{z\in Q}S(z)$ is a zonotope. If the original
parallelotope $P$ is a zonotope, then $P_0$ is a point and $P=Z(Q)$.

Recall that a set of vectors $U$ is a {\em unimodular system} if every
vector $u\in U$ has an integer representation in any basic subset
$B\subseteq U$ (details see in \cite{DaG}). We say that the set of
vectors $Q$ {\em spans} $U$ if there are scalars $\beta_z$, $z\in Q$, such
that $U=\{\beta_zz:z\in Q\}$. It is proved in \cite{Gr}, that the set of
vectors $Q(P)$ of all directions of non-zero width of a parallelotope
$P$ spans a unimodular system $U(P)$. McMullen \cite{McM} proved that a
zonotope $Z(Q)=\sum_{z\in Q}S(z)$ is a parallelotope if and only if $Q$
spans a unimodular system.

We obtain the following result.
\begin{theor}
\label{PZS}
Any parallelotope is:

(i) either a zonotope,

(ii) or a parallelotope of zero width in any direction,

(iii) or the Minkowski sum of a zonotope with a parallelotope of
zero width in any direction.
\end{theor}

\vspace{2mm}
The Proposition~\ref{Pnz} shows that, sometimes, one can add a segment
to a parallelotope $P_0$, in order to obtain another parallelotope. If the
parallelotope $P_0$ has zero width in direction $z$, then obtained
parallelotope $P_0+S(z)$ has another L-type than the original one. In
\cite{Gr} necessary and sufficient conditions are given, when the sum
of a parallelotope with a segment is a parallelotope. For to formulate
these conditions, we introduce some new notions.

Venkov \cite{Ve} proved that a polytope $P$ is a parallelotope if and
only if $P$ itself and all its facets are centrally symmetric and the
projection of $P$ along any $(n-2)$-dimensional face is either
parallelogram, or a centrally symmetric hexagon. The four or six facets,
which are projected into edges of a parallelogram or of a hexagon form
a {\em 2-belt} or a {\em 3-belt}, respectively.

A facet $F$ of a parallelotope $P$ is defined by a {\em facet vector}
$p$, such that the facet $F$ lies in the affine hyperplane
$\{x\in{\bf R}^n:p^Tx=\frac{1}{2}p^Tt\}$. Here $t$ is the {\em lattice
vector} connecting the center of $P=P(0)$ with the center of the
parallelotope $P(t)$ adjacent to $P$ by the facet $F$. Let $I$ be the
set of indices of all pairs of opposite facets of $P$. Then we have:
\begin{equation}
\label{par}
P(0)=\{x\in{\bf R}^n:-\frac{1}{2}p_i^Tt_i\le p_i^Tx\le\frac{1}{2}p_i^Tt_i,
\mbox{  }i\in I\}.
\end{equation}
If $P=P(0)$ is a Voronoi polytope, then the facet vectors $p_i$ are
parallel to the lattice vectors $t_i$, and we can set $p_i=t_i$.
We use the same names 2- and 3-belts for the two and three facet vectors,
defining facets of a 2- and a 3-belt, respectively.

It is proved in \cite{Gr} the following:
\begin{prop}
\label{sum}
For a parallelotope $P$ and a vector $z$, the following assertion are
equivalent:

(i) the Minkowski sum $P+S(z)$ is a parallelotope;

(ii) vector $z$ is orthogonal to at least one facet vector of each
3-belt of $P$.
\end{prop}

\section{The Voronoi polytope of the lattice $D_n$}
We apply Proposition~\ref{sum} to the Voronoi polytope $P_V(D_n)$
of the root lattice $D_n$. The facet vectors of the Voronoi polytope
$P_V(D_n)$ are $n(n-1)$ roots of the root system $D_n$. We take the
roots in the usual form $e_i\pm e_j$, $1\le i<j\le n$. Here
$\{e_i:i\in N\}$, $N=\{1,2,...,n\}$, is an orthonormal basis of
${\bf R}^n$. According to (\ref{par}), we have:
\[P_V(D_n)=\{x\in{\bf R}^n: -1\le x_i\pm x_j\le 1, \mbox{  }
      1\le i<j\le n\}. \]
The vertices of $P_V(D_n)$ are of the following two forms (cf. \cite{CS}):
\[\pm v_i=\pm e_i, i\in N, \mbox{ and }
v(S)= \frac{1}{2}(e(S)-e({\overline S})), S\subseteq N, \]
where ${\overline S}=N-S$ and $e(T)=\sum_{i\in T}e_i$ for any
$T\subseteq N$. The $2^n$ vertices of the set $\{v(S):S\subseteq N\}$
are vertices of a unit cube with its center in origin. The $2n$ vertices
$\pm v_i$, $i\in N$, form $2n$ pyramids having the $2n$ facets of the
unit cube as bases. Hence, the vertex $v_i$ is adjacent to a vertex
$v(S)$ only if $S\ni i$. Similarly, the vertex $-v_i=-e_i$ is adjacent
to a vertex $v(S)$ only if $i\not\in S$.
The edge $v(S)-v_i=\frac{1}{2}e(T)-e({\overline T})$, $T=S-\{i\}$,
connects these vertices. The vertex $v(S)$ is adjacent to a vertex
$v(S')$ only if $v(S)-v(S')=\pm e_i$ for some $i\in N$.

So, edges of $P_V(D_n)$ are of the following form
\begin{equation}
\label{edge}
\mbox{ either }e(S)-e({\overline S}), S\subseteq N,\mbox{ or }e_i, i\in N.
\end{equation}
Up to sign, there are $2^{n-1}+n$ directions of edges, i.e., edge zones.

Each facet of $P_V(D_n)$ is an $(n-1)$-dimensional bipyramid with an
$(n-2)$-dimensional cube as its base. Hence, each 2-face of $P_V(D_n)$ is
a triangle. This implies that all edge zones of $P_V(D_n)$ are open.
\begin{prop}
\label{SDn}
The following assertions are equivalent:

(i) $P_V(D_n)+S(z)$ is a parallelotope;

(ii) $z$ is parallel to an edge of $P_V(D_n)$.
\end{prop}
{\bf Proof}. (i)$\Leftrightarrow$(ii). By Proposition~\ref{sum},
we have to show that the set of vectors $z$, each of which is orthogonal
to at least one vector $p$ of each 3-belt of $P_V(D_n)$ coincides with
the set of edges of $P_V(D_n)$.

The 3-belts of $P_V(D_n)$ are of the following two types:

\[ (a) \mbox{  }e_i-e_j, e_j-e_k, e_i-e_k;\mbox{  }
(b) \mbox{  }e_i+e_j, e_j+e_k, e_i-e_k. \]
We find all vectors $z=\sum_{i=1}^nz_ie_i$ such that $z^Tp=0$ for at
least one facet vector $p$ of each belt. The vector $z$ cannot have 3
mutually non-equal coordinates. In fact, if there are three such
coordinates $z_i$, $z_j$, $z_k$, then $z$ is not orthogonal to any
facet vector of the belt $(e_i-e_j, e_j-e_k, e_i-e_k)$ of type (a).
Therefore, the vector $z$ should be of the form
$z=z'(S):=z_1e(S)+z_2e({\overline S})$, $S\subseteq N$. Since at least
one pair of each triple of indices $\{i,j,k\}$ lies either in $S$, or
in $\overline S$, the vector $z'(S)$ is orthogonal to at least one
vector of each belt of type (a).

For the vector $z'(S)$ to be orthogonal to at least one vector of
each belt of type (b), it should be either $z_1=z_2$, or $z_1=0$ and
$|{\overline S}|=1$, or $z_2=0$ and $|S|=1$, or $z_1+z_2=0$.
The last condition is necessary for the vector $z'(S)$ with $z_1\not=z_2$
to be orthogonal to at least one vector of the belt type (b), such that
either $i\in S$, $k\in\overline S$, or $k\in S$, $i\in\overline S$.

So, we obtain that, up to a multiple, the vector $z$ has the form
\[\mbox{ either } z=z(S)=e(S)-e({\overline S}), S\subseteq N, \mbox{ or }
 z=e_i, i\in N. \]
Comparing these vectors with edges (\ref{edge}) of $P_V(D_n)$, we obtain that all vectors are directed along edges of $P_V(D_n)$. \hfill $\Box$

\section{Four-dimensional unimodular systems and zo\-no\-to\-pal
parallelotopes}
Let $\Sigma_n$ be an $n$-dimensional simplex. It has $\frac{1}{2}n(n+1)$
edges. The $\frac{1}{2}n(n+1)$ vectors which are parallel to edges of
$\Sigma_n$ and have the same length as corresponding edges, form a
maximal unimodular system $A_n$. It represents the graphic matroid of
the complete graph $K_{n+1}$ which is the one-dimensional skeleton of
the simplex $\Sigma_n$. Recall that a set of vectors corresponding to
edges of a graph $G$ represents a graphic (cyclic) matroid of the graph
$G$ if the sum of vectors (taken in suitable direction) along any cycle
of $G$ is zero vector. Changing in this definition cycle by cocycle (cut),
we obtain a representation of the cographic matroid of $G$. (See any book
on Matroid Theory, for example, \cite{Aig}.)

We identify the vectors of $A_n$ with the corresponding edges of
$\Sigma_n$ and $K_{n+1}$. Take the $n$ vectors incident to a vertex
$v\in\Sigma_n$ as a basis of $A_n$ and denote them $e_i$,
$1\le i\le n$. Suppose that these vectors are directed from the vertex
$v$. Then the other $\frac{1}{2}n(n-1)$ vectors of $A_n$ are
$e_i-e_j$, $1\le i<j\le n$.

The Minkowski sum of all vectors of $A_n$ is an $n$-dimensional
zonotope which is called {\em permutohedron}. It is a primitive
parallelotope. Voronoi called its L-type as the {\em principal} type.

Unimodular $n$-dimensional subsystems of $A_n$ are related to
subgraphs of rank $n$ of $K_{n+1}$. Recall that {\em rank} of a graph
is the number of its vertices minus the number of its components.

If a graph $G$ is planar, then there exists its dual planar graph $G^+$
edges of which are in one-to-one correspondence with edges of $G$. The
graphic matroid of $G$ is isomorphic to the cographic matroid of $G^+$.
Both these matroids are represented by a common unimodular system.

A deletion of an element from the cographic matroid of a graph $G$
provides the contraction of the corresponding edge of $G$. It means
that the end vertices of the contracted edge are identified and the
obtained loop is deleted. For the graph $K_{3,3}$, the contraction
of an edge gives a planar graph on 5 vertices. This graph is the
subgraph of $K_5$ obtained by deletion from $K_5$ of two non-adjacent
edges. It is denoted by $K_5-2\times{\bf 1}$. The dual
$(K_5-2\times{\bf 1})^+$ is isomorphic to $K_5-2\times{\bf 1}$.

Here and below instead of the sum ${\bf 1}+{\bf 1}+...+{\bf 1}$ of
$k$ ones of \cite{Co} (denoting $k$ non-adjacent edges), we write
$k\times{\bf 1}$. We use here and below the following Conway's
notations: $C_{ijk...}$ denotes the graph consisting of more than two
chains each containing $i,j,k,...$ edges and all chains connect the same
two vertices. But note that $C_k$ is a cycle with $k$ edges, i.e.
$C_k=C_{ij}$ with $i+j=k$. $G+k\times\bf 1$ denotes a graph $G$ with $k$
pendant edges. For the matroid of the graph $G+k\times\bf 1$ it is not
important, whether the $k$ edges are connected to $G$ or not, or the $k$
edges form a tree or they are disconnected. It is important, that the
subgraph induced by these $k$ edges contains no cycle.

In dimension 4, there are two maximal unimodular systems:

1) $A_4$, representing the graphic matroid of the complete graph
$K_5$, and

2) the unimodular system, representing the cographic matroid $K^*_{3,3}$
of the complete bipartite graph $K_{3,3}$.

There are 16 subgraphs of rank 4 in $K_5$. They are drawn on p.87 of
\cite{Co}. But the graph $K_4$ on this picture has rank 3. It should be
changed by the subgraph $C_{221}+{\bf 1}$ missed in \cite{Co}. A correct
picture of these graphs is given on p.55 of \cite{Va}, but the graph
$C_{221}+{\bf 1}$ is denoted there as $K_4$.

Since proper cographic submatroids of $K^*_{3,3}$ are isomorphic to
graphic ones, in dimension 4, there is only one cographic unimodular
system $K_{3,3}^*$, which is not isomorphic to graphic one. Hence,
besides the mentioned above 16 graphic four-dimensional unimodular
systems, there is the 17th cographic four-dimensional unimodular system
$K^*_{3,3}$. This implies, there are exactly 17 four-dimensional
zonotopal parallelotopes. Amongst them only permutohedron is primitive.

Note that all edge zones of a zonotope are closed. All edges of an edge
zone have the same length. A deletion of a vector from the unimodular
system of a zonotope relates to contraction of the corresponding edge
zone.

The correspondence of zonotopal parallelotopes from \cite{De} with
subgraphs of $K_5$ is given in Table 1. In this table, $N_D$ denotes
the number given to a parallelotope in \cite{De} (we call it
{\em Delaunay number}), and $m$ is the number of segments in the
Minkowski sum of the corresponding zonotope. According to \cite{Co},
${\bf 2}$ and ${\bf 3}$ denote the subgraphs of $K_5$,
which are connected chains of two and three edges, respectively.

{\bf Table 1. Four-dimensional zonotopal parallelotopes $Z(G)$}
\[\begin{array}{|l||l||l|l||l|l||l|l|l|l||} \hline
N_D& 1& 4& 19& 5& 6& 7& 8& 9& 10 \\ \hline
m & 10& 9& 9& 8& 8& 7& 7& 7& 7 \\ \hline
G & K_5 & K_5-{\bf 1}& K^*_{3,3} &K_5-2\times{\bf 1}& K_5-{\bf 2}&
K_5-{\bf 1}-{\bf 2}&K_4+{\bf 1}& C_{2221}&K_5-{\bf 3}\\ \hline
\end{array}\]
\[\begin{array}{|l||l|l|l|l||l|l|l||l||} \hline
N_D& 11& 12& 13&16& 14& 15& 17& 18 \\ \hline
m & 6& 6& 6& 6& 5& 5& 5& 4 \\ \hline
G & C_{222} & C_{321}& C_{221}+{\bf 1} &C_3+C_3& C_4+{\bf 1}&
C_5& C_3+2\times{\bf 1}& 4\times{\bf 1}\\ \hline
\end{array}\]

\section{Unextendible unimodular subsystems of $D_4$}
The four-dimensional Voronoi polytope $P_V(D_4)$ is the self-dual
regular four-dimensional polytope called the {\em 24-cell}. Coxeter
\cite{Cox} denotes it as $\{3,4,3\}$. The edges of $P_V(D_4)$ are given
in (\ref{edge}) where $N=\{1,2,3,4\}$. In this case, up to sign, we have
12 vectors, and these 12 vectors up to the multiple $\sqrt{2}$ are 12
vectors of the root system $D_4$. For convenience, we use below
for edges of $P_V(D_4)$ the usual form of the root system
$D_4=\{e_i\pm e_j: 1\le i<j\le 4\}$. Besides, we denote the vector
$e_i\pm e_j$ by the symbol $ij^{\pm}$, where the signs agree.

Note that $D_4$ consists of the following three 4-sets of mutually
orthogonal vectors: $\{ij^{\pm},kl^{\pm}\}$, where
$\{i,j,k,l\}=\{1,2,3,4\}$. Call such a 4-set {\em quadruple}. Each
quadruple relates to one of the three partitions of the 4-set
$\{1,2,3,4\}$ into pairs. Call a three mutually orthogonal roots of $D_4$
by a {\em triple}. Each triple $t$ is contained in a uniquely determined
by $t$ quadruple $q_t\supset t$.

Consider the following two maximal unimodular systems contained in the
root system $D_4$: the graphic system $A_4-e$, where $e$ is any vector of
$A_4$, and the cographic system $K_{3,3}^*$. The system $A_4-e$ represents
the graphic matroid of the graph $K_5-{\bf 1}$, and the system $K^*_{3,3}$
represents the cographic matroid of the complete bipartite graph $K_{3,3}$.
There are many ways to choose in $D_4$ vectors forming the unimodular
systems $A_4-e$ and $K^*_{3,3}$. We choose the vectors as follows.
The 6 vectors $ij^-$, $1\le i<j\le 4$, form the graphic system $A_3$
representing the graphic matroid of the complete graph $K_4$. If the
vertices of $K_4$ are denoted by the numbers 1,2,3,4, then the vector
$ij^-=e_i-e_j$ represents the edge $(ij)$ connecting the vertices $i$
and $j$. Now suppose that the vertex 5 of $K_5-{\bf 1}$ is not connected
with the vertex 1 of its subgraph $K_4$. Then we can relate the vectors
$12^+$, $13^+$, $14^+$ to the edges (25), (35), (45), respectively. It
is easy to verify that the 9 vectors $ij^-$, $1\le i<j\le 4$, $1i^+$,
$2\le i\le 4$, form the unimodular system $A_4-e$. This unimodular
system consists of the following three triples of mutually orthogonal
vectors: $(ij^-,ij^+,kl^-)$, $ij=12,13,14$, corresponding to three
partitions of the 4-set $\{i,j,k,l\}$ into pairs.

The graph $K_5-{\bf 1}$ is planar. It has three vertices of degree 4
and two vertices of degree 3. The 9 edges of this graph are partitioned
into the following two orbits (1) and (2) of the automorphism group of
$K_5-{\bf 1}$:

(1) 3 edges with both end vertices of degree 4; they are represented by
the roots $23^-$, $24^-$, $34^-$;

(2) 6 edges with end vertices of degree 3 and 4; they are represented by
the roots $12^{\pm}$, $13^{\pm}$, $14^{\pm}$.

If we delete in the graph $K_5-{\bf 1}$ the edge (24), we obtain the
planar graph $G_5:=K_5-2\times{\bf 1}$. The vertices 1, 2, 4, 5 of the
graph $G_5$ have degree 3 and form a 4-cycle with edges (12), (25),
(45), (14). The vertex 3 has degree 4 and it is adjacent to the four
vertices of the 4-cycle by edges (13), (23), (35), (34). Now, we
consider the dual graph $G_5^+$ which is isomorphic to original one.
Let the vertices of the 4-cycle of the $G_5^+$ be $a, b, c, d$ in this
order along the 4-cycle. Then edges of this 4-cycle, corresponding to
the edges (13), (23), (35), (34) of $G_5$, have the following pairs of
end vertices, respectively: $(ab), (bc), (cd), (ad)$. The vertex $v$ of
degree 4 is adjacent to the vertices $a, b, c, d$ by edges, corresponding
to the edges (14), (12), (25), (45) of $G_5$, respectively.
The four last edges form a 4-cocycle of $G_5^+$. The 8 vectors $12^-$,
$13^-$, $14^-$, $23^-$, $34^-$, $12^+$, $13^+$, $14^+$ related to the
edges of $G_5$ and $G_5^+$ represent the graphic matroid of $G_5$ and the
cographic matroid of $G_5^+$. The cographic matroid of $G_5^+$ is a
submatroid of $K^*_{3,3}$.

The graph $G_5^+$ is obtained from $K_{3,3}$ by the contraction of one
of its edges. The operation opposite to the contraction of an edge of
the graph $K_{3,3}$ is the splitting of the vertex $v$ of degree 4 in
$G_5^+$ into two adjacent vertices $v'$ and $v''$. This splitting is such
that the vertex $v'$ is adjacent to the vertices $a$ and $c$, and the
vertex $v''$ is adjacent to the vertices $b$ and $d$. We obtain the
complete bipartite graph $K_{3,3}$ with mutually non-adjacent vertices
$b,d,v'$ of one part and mutually non-adjacent vertices $a,c,v''$ of
other part. The vectors related to edges $(av')$, $(cv')$ and $(bv'')$,
$(dv'')$ are $14^-$, $12^+$ and $12^-$, $14^+$, respectively. Since
$12^+-14^-=14^+-12^-=24^+$, we should relate the vector $24^+$ to the
edge $(v'v'')$.

We see that this representation of the cographic matroid $K^*_{3,3}$
is obtained from the representation of the graphic matroid of the
graph $K_5-{\bf 1}$ by changing the vector $24^-$ into the vector
$24^+$.

Call a unimodular subsystem $U\subseteq D_4$ {\em unextendible} if
the set $U\cup\{r\}$ is not unimodular for all $r\in D_4-U$.
\begin{prop}
\label{3uni}
Up to isomorphism, the root system $D_4$ contains exactly the following
3 unextendible unimodular subsystems: a quadruple, $A_4-e$, and
$K_{3,3}^*$.
\end{prop}
{\bf Proof}. Let $q=\{ij^-,ij^+,kl^-,kl^+\}\subset D_4$ be a quadruple.
Obviously it is a unimodular system. It is a basis of the space
${\bf R}^4$. Any root $r\in D_4-q$ has a unique non-integer representation
in this basis. So, the system $q\cup\{r\}$ is not unimodular for all
$r\in D_4-q$. Hence, each quadruple of $D_4$ is an unextendible unimodular
subsystem of $D_4$.

This implies that any other unimodular subsystem of $D_4$ does not
contain a quadruple. Consider a set $U=t_1\cup t_2\cup t_3$ of three
mutually disjoint triples. Obviously, $U$ is a maximal subset of $D_4$
not containing a quadruple. (Note, there are $4^3=64$ such sets.)
We show that $U$ is a unimodular system isomorphic either to $A_4-e$, or
to $K_{3,3}^*$. Obviously, these systems are unextendible.

Each triple $t\subset U$ is complemented by a unique root $r_t\in D_4$
up to the quadruple $q_t$. Label the set $U$ by the {\em triad}
$(r_1,r_2,r_3)$ of these complementing roots. For example, the
following set $U=\{ij^-,ij^+,kl^-:ij=12,13,14\}$ is labeled by the
triad $(34^+,24^+,23^+)$.

The automorphism group of the root system $D_4$ consists of the following
operations: permutations of indices, changing a pair of indices $(ij)$
by the complementing pair $(kl)$, reversing the sign of a unit vector
$e_i\to -e_i$. The collection of all triads is partitioned into two
orbits of the automorphism group. One orbit consists of triads with even
number of minus signs. Another orbit contains triads with odd number of
minus signs. The sets $U$ labeled by triads of these two orbits are
isomorphic to unimodular systems $A_4-e$ and $K_{3,3}^*$, respectively.
\hfill $\Box$

Denote the zonotopes related to a unimodular system $U$ by $Z(U)$. Since
the deletion of the vector $24^-$ from $A_4-e$ and the vector $24^+$ from
$K_{3,3}^*$ provides the same unimodular system $U(K_5-2\times{\bf 1})$,
the contraction of the edge zones of $Z(A_4-e)$ and $Z(K^*_{3,3})$
corresponding to $24^-$ and $24^+$, respectively, provides zonotopes,
which are both isomorphic to $Z(U(K_5-2\times{\bf 1}))$.

\section{Sums of $P_V(D_4)$ with zonotopes}
The sums $P_V(D_4)+Z(A_4-e)$ and $P_V(D_4)+Z(K^*_{3,3})$ are
primitive parallelotopes. Their projections along a closed edge zone
into three-dimensional space are drawn in Figs.II and III of \cite{De};
see also \cite{De1}.
Their Delaunay numbers are $N_D=2$ and $N_D=3$. In Figs.II and III,
the closed edge zones of these parallelotopes are denoted by numbers
1,2,3,...,8,9. The number 9 corresponds to the edge zone, along which
the parallelotope is projected. Any non-primitive parallelotope is
obtained from these two ones by contracting some closed edge zones.
Some edges are denoted by $0i$, $1\le i\le 8$. This means that after
the contraction of the edge zone with number $i$ the edge with number
$0i$ is contracted to an edge of $P_V(D_4)$ which is denoted by 0.

Figs.II and III show that the contraction of the edge zone 4 in
$P_V(D_4)+Z(A_4-e)$ and the edge zone 6 in $P_V(D_4)+Z(K^*_{3,3})$
gives the same parallelotope. Comparing the numbers of edge zones of this
parallelotope in Figs.II and III, we obtain their correspondence, shown
in Table below. Besides, in this Table the vectors of the unimodular
systems representing $A_4-e$ and $K^*_{3,3}$ are given.

\[\begin{array}{|c||c|c|c|c|c|c|c|c|c|c|} \hline
{\rm Fig.I}&1&2&3&4&5&6&7&8&9&-\\  \hline
{\rm Fig.II}&2&1&5&-&7&8&3&4&9&6\\ \hline
{\rm roots}&23^-&14^-&12^-&24^-&34^-&13^-&12^+&13^+&14^+&24^+\\ \hline
\end{array}\]

We see that the edge zones 1, 4 and 5 of Fig.I contains edges of the
same orbit (1) of the automorphism group of $K_5-{\bf 1}$. Hence the
contraction in $P_V(D_4)+Z(A_4-e)$ of any edge zone of this orbit gives
isomorphic parallelotopes. The parallelotope with contracted edge zone 1
is the parallelotope of \cite{De} with the Delaunay number $N_D=21$.
This parallelotope is $P_V(D_4)+Z(U(K_5-2\times{\bf 1}))$.

Note that for zonotopes $Z(U)$, which are parallelotopes, it is not
important whether the summing vectors are orthogonal or not. A
parallelepiped and a cube have the same L-type. But an orthogonality of
summing vectors in $Z(U)$ affects heavenly onto the L-type of the sum
$P_V(D_4)+Z(U)$. This implied by the fact that any two orthogonal roots
do not belong to a 2-face of $P_V(D_4)$. Hence the sum
$P_V(D_4)+S(r)+s(r')$, where $r$ and $r'$ are orthogonal, obtains a new
2-face spanned by $r$ and $r'$ contrary to the case, when $r$ and $r'$
are not orthogonal.

Facets of $P_V(D_4)$ are octahedra whose edges are roots of $D_4$.
Each octahedron has four pairs of opposite parallel triangle faces and
six pairs of parallel edges, which are parallel to six distinct roots.
This six roots are partitioned into three pairs of orthogonal roots. Each
triangle contains one representative root from these three pairs. Each
of the six roots belongs to four triangles.

For a facet $F$ of $P_V(D_4)$, let $R(F)$ be the set of the six roots
parallel to edges of $F$. In this case, when edges are parallel to roots,
the facet vectors are given by (\ref{edge}). The facet vector $e_k$,
$1\le k\le 4$, defines a facet $F$, such that
$R(F)=\{ij^{\pm}:i,j\not=k\}$.
This is the root system $D_3$, which is isomorphic to $A_3$. The facet
vector of another type $\frac{1}{2}\sum_{i=1}^4\varepsilon_ie_i$, where
$\varepsilon_i\in\{\pm 1\}$, defines a facet $F$, such that
$R(F)=\{ij^{-\varepsilon_i\varepsilon_j}:1\le i<j\le 4\}$.

Note that, for each pair $(r,r')$ of orthogonal roots, there is a facet
$F$ of $P_V(D_4)$, such that $r,r'\in R(F)$.

Let $U=U(G)$ be a unimodular system of roots representing a subgraph
$G\subseteq K_5-\bf 1$ or $U=K_{3,3}^*$. Let $\pi(U)$ be the set of
maximal pairs of orthogonal roots in $U$. A pair $(r,r')\subseteq U$
of orthogonal roots in $U$ is called {\em maximal} if there is no
root in $U-\{r,r'\}$ which is orthogonal to both the roots $r,r'$.
Let $\tau(U)$ be the set of triples of mutually orthogonal roots in $U$.
Recall that each triple $t\in \tau(U)$ uniquely determines the fourth
root $r_t\in D_4$, such that $r_t$ is orthogonal to all roots of $t$.

Recall that $P_V(D_4)$ has 24 facets (it is a 24-cell) and all its
16 belts are 3-belts.
\begin{prop}
\label{PVZ}
Let $Z(U)=\sum_{r\in U}S(r)$. Then it holds:

(i) the sum $P_V(D_4)+Z(U)$ has $|\pi(U)|$ 2-belts and $16+3|\tau(U)|$
3-belts;

(ii) the sum $P_V(D_4)+Z(U)$ has $24+2|\tau(U)|$ facets;

(iii) the sum $P_V(D_4)+Z(U)+S(r)$ is a parallelotope for all $r\in D_4$
such that $r\not=r_t$, $t\in \tau(U)$.
\end{prop}
{\bf Proof}. Let $F$ be a facet of $P_V(D_4)$ and $r\in R(F)$. The facet
$F$ is transformed to the facet $F+S(r)$ in the sum $P_V(D_4)+S(r)$.
(Here $S(r)$ is a segment of the line parallel to the root $r$.) The facet
$F+S(r)$ has also 4 pairs of parallel faces. But the four triangles of $F$
containing an edge parallel to $r$ are transformed into trapezoids with
two edges parallel to $r$.

Now consider the sum $P(r,r'):=P_V(D_4)+S(r)+S(r')$, $r,r'\in R(F)$.
If $r$ and $r'$ are not orthogonal, then the facet $F+S(r)+S(r')$ has
also 4 pairs of parallel faces. The two parallel faces which contains
edges parallel to both the roots $r$ and $r'$ are transformed into
pentagons. The two pairs of faces having edges parallel only one of
these two roots are transformed into two pairs of trapezoids. One pair of
parallel faces is not changed.

If $r$ and $r'$ are orthogonal, then each face of $F$ has only one edge
parallel to one of these roots. Hence, each face is transformed into a
trapezoid. There are two opposite vertices of $F$ which are not
incident to the edges parallel to these roots. These vertices are
transformed into new square faces $Q$ of $F+S(r)+S(r')$. So, this facet
has now 10 faces.

This shows that the L-type of the sum $P_V(D_4)+S(r)+S(r')$
depends on whether $r$ and $r'$ are orthogonal or not. Besides, if
$r$ and $r'$ are orthogonal, then there are two pairs of opposite
facets which are transformed into polyhedra with 10 faces. These four
facets form a new 2-belt. So, each pair of orthogonal roots in $U$
generates a 2-belt.

Let $t=(r,r',r'')\in \tau(U)$. Consider the sum $P_t:=P(r,r')+S(r'')$.
Let $F$ be a facet such that $r,r'\in R(F)$. Then $r''\not\in R(F)$.
Let $Q$ be the quadrangle face of $F+S(r)+S(r')$. Then $Q+S(r'')$ is a
cube. It is a new facet of $P_t$. This facet, its opposite and the
facets of the 2-belt of $P(r,r')$ form a new 3-belt $\cal B$. But the
cube $Q+S(r'')$ has 3 pairs of opposite faces, and each pair of its
faces generates a 3-belt. Hence, for each $t\in \tau(U)$, the parallelotope
$P_t$ has additionally three belts. We proved (i) and (ii).

For $t\in \tau(U)$, the root $r_t$ is orthogonal to the new facet
$Q+S(r'')$. So, $r_t$ is the facet vector of this facet. Obviously,
$r_t$ is not orthogonal to any facet vector of the new 3-belt $\cal B$.
By Proposition~\ref{sum}, $P_t+S(r_t)$ is not a parallelotope.
Since the facet vector of the transformed facet $F+Z(U)$ is the same
as the facet vector of $F$, we obtain (iii). \hfill $\Box$

\vspace{2mm}
Table 2 shows the 35 unimodular systems $U$ of roots sums of which
with $P_V(D_4)$ are parallelotopes. Maximal pairs and triples of mutually
orthogonal roots are distinguished by parentheses. Besides of used above
denotations of graphs, $H_k$ denotes the skeleton of a $k$-dimensional
parallelepiped. As in Table 1, $N_D$ denotes the Delaunay number of the
parallelotope $P_V(D_4)+Z(U)$ and $m$ is the number of roots in $U$.
The parallelotope, missed by Delaunay and found by Shtogrin, is denoted by
{\bf St}. The fifth column gives dimension of the added zonotope $Z(U)$.
Note that dim$Z(U)$ is equal to the rank of the graph $G$, which is
represented by the unimodular system $U$.

For to compare Tables 2 and 1, we add the last column. In this column
$N_D^0$ denotes the Delaunay number $N_D$ of the corresponding zonotopal
parallelotope $Z(U)$ if dim$U$=4; $N_D^0=a_i$, $1\le i\le 5$, denotes
a 3-dimensional zonotopal parallelotope if dim$U$=3; and
$N_D^0=\alpha, \beta_1, \beta_2$ denotes a 2-dimensional zonotopal
parallelotope if dim$U$=2. Note that

$a_1$ denotes a permutohedron = a truncated octahedron;

$a_2$ denotes an elongated dodecahedron;

$a_3$ denotes a prism with a hexagonal base;

$a_4$ denotes a rhombic dodecahedron;

$a_5$ denotes a parallelopiped with mutually orthogonal edges;

$a'_5$ denotes a parallelopiped with a pair of parallel rectangle facets;

$a''_5$ denotes a parallelopiped without rectangle facets;

$\alpha$ denotes a centrally symmetric hexagon;

$\beta_1$ denotes a rectangle;

$\beta_2$ denotes a parallelogram without orthogonal edges.

\newpage
{\bf Table 2. Zonotopes $Z(U)$ such that $P_V(D_4)+Z(U)$ is a parallelotope}
\[\begin{array}{|l|l|l|l|c||l|} \hline
N_D& m &  \mbox{roots of the unimodular system }U &\mbox{graph}
& \mbox{dim}U & N_D^0 \\ \hline
- & 10 &  A_4 & K_5 & 4 & 1 \\ \hline
2  & 9 & (12^-,12^+,34^-),(13^-,13^+,24^-),(14^-,14^+,23^-)& K_5-{\bf 1}
& 4 & 4\\
3  & 9 & (12^-,12^+,34^-),(13^-,13^+,24^+),(14^-,14^+,23^-)& K_{3,3}^*
& 4 & 19\\
\hline
20 & 8 & (12^-,12^+,34^-),(13^-,13^+,24^-),(14^+,23^-)& K_5-{\bf 2}& 4
& 6\\
21 & 8 & (12^-,12^+,34^-),(13^-,13^+,24^-),(14^-,14^+)&
K_5-2\times{\bf 1}& 4 & 5\\ \hline
22 & 7 & (12^+,34^-),(13^-,13^+,24^-),(14^-,14^+)& K_5-{\bf 3}& 4&10 \\
23 & 7 & 34^-,(13^-,13^+,24^-),(14^-,14^+,23^-)& C_{2221}& 4& 9 \\
24 & 7 & (12^-,12^+,34^-),(13^-,13^+,24^-),14^+& K_5-{\bf 1}-{\bf 2}& 4
& 7\\
25 & 7 & (12^+,34^-),(13^-,13^+,24^-),(14^+,23^-)& K_4+{\bf 1}& 4& 8 \\
26 & 7 & (12^-,12^+,34^-),(13^-,13^+),(14^-,14^+)& K_5-{\bf 1}-{\bf 2}&
4 & 7\\ \hline
27 & 6 & (12^-,12^+,34^-),(13^-,13^+),14^+ &  C_{321}& 4 & 12\\
28 & 6 & (12^+,34^-),(13^-,13^+,24^-),14^+ & C_{221}+{\bf 1}& 4 &13 \\
29 & 6 & (13^-,13^+,24^-),(14^-,14^+,23^-) & C_{222}& 4 & 11\\
30 & 6 & (12^+,34^-),(13^+,24^-),(14^-,14^+)&C_{221}+{\bf 1}& 4 & 13\\
31 & 6 & (12^-,12^+),(13^-,13^+),(14^-,14^+)& C_{222} & 4 & 11 \\
32 & 6 & (12^+,34^-),(13^-,24^-),(14^-,14^+)& C_3+C_3 & 4 & 16\\
33 & 6 & (12^+,34^-),(13^+,24^-),(14^+,23^-) & K_4 & 3 & a_1\\ \hline
34 & 5 & 12^+,(13^-,13^+,24^-),14^+ & C_3+2\times{\bf 1} & 4 & 17\\
35 & 5 & (12^-,12^+,34^-),13^-,14^+ & C_5 & 4 & 15 \\
36 & 5 & (13^-,13^+,24^-),(14^-,14^+) & C_4+{\bf 1} & 4 & 14\\
37 & 5 & (12^-,12^+),(13^-,13^+),14^+ & C_4+{\bf 1} & 4 & 14\\
38 & 5 & (12^+,34^-),(13^-,13^+),14^+ & C_3+2\times{\bf 1} & 4 &17 \\
39 & 5 & (12^+,34^-),(13^+,24^-),14^+ & C_{221} & 3 & a_2\\ \hline
40 & 4 & (13^-,13^+,24^-),14^+& 4\times{\bf 1}=H_4 & 4 & 18 \\
41 & 4 & 12^+,(13^-,13^+),14^+& 4\times{\bf 1}=H_4 & 4 & 18 \\
42 & 4 & (12^+,34^-),13^+,14^+& C_3+{\bf 1} & 3 & a_3\\
43 & 4 & (13^+,24^-),(14-,14^+)& 4\times{\bf 1}=H_4 & 4 & 18 \\
44 & 4 & (13^-,13^+),(14^-,14^+)& C_4 & 3 & a_4\\ \hline
45 & 3 & (14^-,14^+,23^-)&3\times{\bf 1}=H_3 & 3 & a_5\\
46 & 3 & (13^-,13^+),14^-& 3\times{\bf 1}=H_3 & 3 & a'_5\\
47 & 3 & 12^+,13^+,14^+& 3\times{\bf 1}=H_3 & 3 & a''_5\\
{\bf St} & 3& 34^-,13^+,14^+& C_3 & 2 & \alpha \\ \hline
48 & 2 & (14^-,14^+)& 2\times{\bf 1}=H_2 & 2 & \beta_1\\
49 & 2 & 13^+,14^+& 2\times{\bf 1}=H_2 & 2 & \beta_2\\ \hline
50 & 1 & 14^+& {\bf 1}=H_1 & 1& \\ \hline \hline
51 & 0 &     & {\rm 24-cell} & &\\ \hline
\end{array} \]

\end{document}